\documentstyle[12pt]{amsart}

\newtheorem{prop}{Proposition}[section]
\newtheorem{dfn}[prop]{Definition}
\newtheorem{theo}[prop]{Theorem}

\newtheorem{rem}[prop]{Remark}
\newtheorem{coro}[prop]{Corollary}

\newtheorem{lem}[prop]{Lemma}

\def\C{{\bold C }}
\def\R{{\bold R }}
\def\Re{ \operatorname{Re} }
\def\Im{ \operatorname{Im} }
\def\Z{{\bold Z }}
\def\Q{{\bold Q }}

\def\a{ \alpha }
\def\b{ \beta }
\def\d{ \delta }
\def\e{ \varepsilon }
\def\g{ \gamma}

\def\l{ \lambda }
\def\L{ \Lambda }

\def\s{ \sigma }
\def\S{ \Sigma }
\def\G{ \Gamma }
\def\p{ \varphi }

\def\m{ \mu }

\def\ra{\rightarrow}
\def\da{\downarrow}
\def\Val{\operatorname{Val}}

\def\dim{\operatorname{dim}}

\title{Height zeta functions of toric varieties}

\author{Victor V. Batyrev}
\thanks{Supported by Deutsche Forschungsgemeinschaft} 
\address{Universit\"at-GHS-Essen, Fachbereich  6,  Mathematik \\
 Universit\"atsstr. 3,  45141  Essen, Germany  }
\email{victor.batyrev@@uni-essen.de }

\author{Yuri Tschinkel}
\thanks{Leibniz Fellow at ENS, Paris}
\address{Dept. of Mathematics, U.I.C.\\
Chicago, (IL) 60608,  U.S.A.}
\email{yuri@@math.uic.edu}

\begin{document}

\date{}

\maketitle

\thispagestyle{empty}

\begin{abstract}
We investigate  analytic properties of height zeta 
functions of toric varieties. Using the height   
zeta functions, we prove an asymptotic formula 
for the number of rational points of bounded height 
with respect to an arbitrary line bundle whose first 
Chern  class is contained in the interior of the cone 
of effective divisors. 
\end{abstract}

\newpage

\tableofcontents

\bigskip
\bigskip
\bigskip

\vskip 0,5cm
\section{Introduction}

\bigskip 

Let $X$ be a $d$-dimensional 
algebraic variety defined over a number field $F$.
Denote by ${\cal L}=(L,\{\|\cdot\|_v\})$ a metrized line bundle on $X$ ,
i.e. a line bundle $L$ together with a family of $v$-adic metrics,
where $v$ runs over the set ${\Val }(F)$ of all valuations of $F$. 
For any locally closed algebraic subset  
$Y\subset X$ we denote by $Y(F)$ the set of $F$-rational points in $Y$.
A metrized line bundle  ${\cal L}$ defines a height function
$$
H_{\cal L} \, :\, X(F) \ra \R_{>0}.
$$ 

Assume that a subset $Y \subset X$ and a bundle $L$ are choosen
in such a way that  
$$
N(Y,{\cal L},B) := \# \{ x\in Y(F) \mid H_{\cal L}(x)\le B \} < \infty
$$
for all $B \in \R_{>0}$ (e.g., this holds for any $Y\subset X$ if $L$ is
ample). Then the asymptotic behavior of 
the counting function $N(Y,{\cal L},B)$ as $B \ra \infty$  
is determined by analytic properties of the 
{\em height zeta function} $Z(Y,{\cal L},s)$ 
defined by the series
$$
Z(Y,{\cal L},s) : =\sum_{x\in Y(F)} H_{\cal L}(x)^{-s}
$$
which converges for ${\Re }(s)\gg 0$.
More precisely, one has the following  Tauberian statement:

\begin{theo} {\rm \cite{BaTschi1}}
Assume that the series $Z(Y,{\cal L},s)$ is absolutely
convergent for ${\Re }(s) > a> 0$ and that there exists 
some positive integer $b$ such that  
$$
Z(Y,{\cal L},s)=\frac{g(s)}{(s-a)^b} + h(s)
$$
where $g(s)$ and $h(s)$ are functions holomorphic in the domain
${\Re }(s)\ge a$ and $g(a)\neq 0$. Then 
the following asymptotic formula holds:
$$
N(Y,{\cal L},B)=\frac{g(a)}{a(b-1)!}
B^{a}(\log B)^{b-1}(1+o(1))\;{ for}
\;B\ra \infty. 
$$
\end{theo}

\begin{dfn}
{\rm Denote by $NS(X)$ be the Neron-Severi group of $X$ and by
$NS(X)_{\R} =
NS(X)\otimes \R$. The {\em cone of effective
divisors of  $X$} is the closed cone 
$\L_{\rm eff}(X)\subset NS(X)_{\R}$ generated by the classes of 
effective divisors. }
\end{dfn}

\noindent
Let $[L] \in NS(X)$ be the first Chern class of $L$. 
Denote by ${\cal K}_X = (K_X,\{\|\cdot\|_v\})$ 
the metrized canonical line
bundle on $X$. 

\begin{dfn}
{\rm Let $L$ be any line bundle on $X$. Define  
$$
a(L) : =\inf \{ a\in \R \,|\, a[L]+[K_X]\in \L_{\rm eff}(X)\}.
$$}
\end{dfn}
One of our main results in this paper is the following theorem:
\begin{theo} Let
$T$ be a $d$-dimensional
algebraic torus over a number field $F$,  
$X$ a smooth projective toric variety containing $T$
as a Zariski open subset, and  
${\cal L}$ a metrized line bundle on $X$ $($with the metrization
introduced in {\rm \cite{BaTschi1}}$)$. Assume that the class $[L]$ is
contained in the interior of the cone of effective divisors. 
Then the height zeta function has the following representation: 
$$
Z(T,{\cal L},s)=\frac{g(s)}{(s-a(L))^{b(L)}} + h(s)
$$
where $g(s)$ and $h(s)$ are functions holomorphic in the domain
${\Re }(s)\ge a(L)$, $g(a(L))\neq 0$,  and $b(L)$ is the codimension of 
the minimal face of $\L_{\rm eff}(X)$ which contains the
class $a(L)[L]+[K_X]$.  
\label{t1}
\end{theo}

\begin{coro} 
We have the following asymptotic formula
$$
N(T,{\cal L},B)=
\frac{g(a(L))}{a(L)(b(L)-1)!}B^{a(L)}(\log B)^{b(L)-1}(1+o(1))\;{ for}
\;B\ra \infty.
$$
\end{coro}
\noindent
The paper is organized as follows: 
\medskip

The  technical heart of the paper is 
contained in Section 2, where we investigate analytic 
properties of some complex valued functions related 
to convex cones. 

In Section 3, we review  basic facts from harmonic analysis on
the adele group of an algebraic torus. 

In Section 4, we recall the terminology from the theory 
of toric varieties as well as the definition 
and main properties of heights on toric varieties. 

In Section 5, we give the proof of \ref{t1}.  We remark that the 
most subtle part  in the statement of \ref{t1} is the 
nonvanishing of the asymptotic constant $g(a(L)) \neq 0$. 
\bigskip

\section{Technical theorems}

Let $I$ and $J$ be two positive integers, 
${\R}\lbrack {\bold s},{\bold t}  \rbrack$
(resp. ${\C}\lbrack {\bold s}, {\bold t} \rbrack$)
the ring of polynomials in $I +J$  
variables $s_1, \ldots, s_I, t_1, \ldots, t_J$
with coefficients in
${\R}$ (resp. in ${\C}$) and  ${\C}\lbrack \lbrack {\bold s}, {\bold t}
\rbrack \rbrack$  the ring of formal power series in $s_1,
\ldots, s_I,t_1, \ldots, t_J$
with complex coefficients.

\begin{dfn}
{\rm Two elements $f({\bold s},{\bold t}),\, g({\bold s},{\bold t}) 
\in {\C}\lbrack \lbrack {\bold s}, {\bold t}
\rbrack \rbrack$ will be called {\em coprime}, if
$g.c.d.(f({\bold s},{\bold t}),\, g({\bold s},{\bold t})) =1$.
}
\end{dfn}

\begin{dfn}
{\rm Let  $f({\bold s},{\bold t})$ be an  element of 
${\C}\lbrack \lbrack {\bold s}, {\bold t}
\rbrack \rbrack$. By the {\em order} of a monomial
$s_1^{\a_1} \cdots s_I^{\a_I}t_1^{\b_1} \cdots t_J^{\b_J}$ 
we mean the sum
of the exponents 
$$\a_1+ \cdots + \a_I + \b_1 + \cdots + \b_J.$$ 
By the  {\em multiplicity $\mu(f({\bold s},{\bold t}))$
of $f({\bold s},{\bold t})$ at
${\bold 0} = (0, \ldots, 0)$} we always mean
the minimal order of non-zero monomials appearing in the
Taylor
expansion of $f({\bold s},{\bold t})$ at ${\bold 0}$ .   }
\label{mult1}
\end{dfn}

\begin{dfn}
{\rm
Let  $f({\bold s},{\bold t})$ be  a  meromorphic at ${\bold 0}$ function.
Define the   {\em multiplicity $\mu(f({\bold s},{\bold t}))$ }
of $f({\bold s},{\bold t})$ at
${\bold 0}$ as
\[  \mu(f({\bold s},{\bold t})) = \mu(g_1({\bold s},{\bold t})) - 
\mu(g_2({\bold s},{\bold t})) \]
where $g_1({\bold s},{\bold t})$ and $g_2({\bold s},{\bold t})$ are two coprime
elements in  ${\C}\lbrack \lbrack {\bold s},{\bold t} \rbrack \rbrack$
such that $f = g_1/g_2$.
}
\label{mult2}
\end{dfn}

\begin{rem}
{\rm It is easy to show that for any two
meromorphic at ${\bold 0}$ functions  $f_1({\bold s},{\bold t})$ and 
$f_2({\bold s},{\bold t})$,
one has

(i) $\mu(f_1 \cdot f_2) = \mu(f_1) + \mu(f_2)$ (in
particular, one can omit "coprime" in Definition \ref{mult2});

(ii) $\mu(f_1 + f_2) \geq \min \{ \mu(f_1),  \mu(f_2) \}$;

(iii)  $\mu(f_1 + f_2) = \mu(f_1)$ if $\mu(f_2) > \mu(f_1)$. }
\label{mult3}
\end{rem}

Using the  properties \ref{mult3}(i)-(ii), one immediately
obtains from  Definition \ref{mult1} the following:

\begin{prop}
Let $f_1({\bold s},{\bold t}) \in {\C}\lbrack \lbrack {\bold s}, {\bold t}
\rbrack \rbrack$ and $f_2({\bold s}) \in 
{\C}\lbrack \lbrack {\bold s}\rbrack \rbrack$ 
be two analytic at ${\bold 0}$ functions, 
$l({\bold s}) \in {\bold R}[{\bold s}]$ 
a homogeneous linear function in the variables $s_1, \ldots, 
s_I$, 
$${\bold \g} = ({\bold \g}^I, {\bold \g}^J)  
= (\g_1, \ldots, \g_I, \g_1', \ldots, \g_J') \in {\C}^{I
+J}$$ 
an arbitrary complex
vector with $l({\bold \g}^I) \neq 0$,
and $g({\bold s},{\bold t}) := f_1({\bold s},{\bold t})/f_2({\bold s})$. Then
the multiplicity  of the function 
$$
 \tilde{g}({\bold s},{\bold t}): = \left(\frac{\partial}{\partial
z}\right)^k
g({\bold s} + z \cdot {\bold \g}^I, {\bold t} + z\cdot  {\bold \g}^J ) 
|_{z = - l({\bold s})/l({\bold \g}^I)}
$$
at ${\bold 0}$ is at least $\mu(g) - k$, if 
$$
f_2({\bold s} + z \cdot \g^I) |_{z = - l({\bold s})/l({\bold \g}^I)}
$$
is not identically zero. 
\label{mult4}
\end{prop}

Let $\Gamma \subset {\Z}^{I+J}$ be a
sublattice, $\Gamma_{\R} \subset {\R}^{I + J}$ (resp.
$\Gamma_{\C} \subset {\C}^{I +J}$) the scalar extension of
$\Gamma$ to a ${\R}$-subspace
(resp. to a ${\C}$-subspace). 
We always assume that $\Gamma_{\R}
\cap {\R}_{\geq 0}^{I +J} = 0$ and $\G_{\R} \cap \R^J = 0$. 
We set  $P_{\R}: = 
{\R}^{I +J}/\Gamma_{\R}$ and
$P_{\C}: = {\C}^{I+J}/\Gamma_{\C}$. Let $\pi^I$ be the 
natural projection $\C^{I+J} \ra \C^I$. Denote by $\psi$ 
(resp. by $\psi^I$) the canonical surjective mapping 
$\C^{I+J} \ra P_{\C}$ (resp. ${\bold C}^{I} \ra \C^I/\pi^I(\G_{\C})$).

\begin{dfn}
{\rm A complex  analytic function 
$$
f({\bold s}, {\bold t})=
f(s_1, \ldots, s_I, t_1, \ldots, t_J): U \ra {\bold C}
$$
defined
on an open subset $U \subset {\bold C}^{I +J}$
is said to {\em descend to $P_{\C}$} if 
$ f({\bold u} ) = f({\bold u}')$ 
for all ${\bold u}, {\bold u}' \in U$ with ${\bold u} - {\bold u}' \in 
\G_{\C}$. }
\end{dfn}

\begin{rem}
{\rm By definition, if $f({\bold s}, {\bold t})$ 
descends to $P_{\C}$, then there exists
an analytic function $g$ on $\psi(U) \subset P_{\C}$ such that
$f = g \circ \psi$. Using Cauchy-Riemann equations, one
immediatelly   obtains that $f$ descends to $P_{\C}$ if and only
if for any vector ${\bold \a} \in \Gamma_{\R}$ and
any ${\bold u}= (u_1, \ldots, u_{I+J}) \in U$ such that 
${\bold u} + i{\bold \a} \in U$, one has
\[ f({\bold u}+ i{\bold \a})  = f({\bold u}).  \]}
\label{desc}
\end{rem}

\begin{dfn}
{\rm An analytic function $f({\bold s}, {\bold t})$ in the domain
${\Re }({\bold s}) \in {\R}_{>0}^I$, ${\Re }({\bold t}) 
\in {\R}_{> -\d_0}^J$ $(${for some} $\d_0 >0)$ is called 
{\em good with respect to $\Gamma$ and
the set of variables $\{s_1,...,s_I\}$}
if it satisfies the following conditions:

{(i)} $f({\bold s},{\bold t})$ descends to $P_{\C}$;

{(ii)} There exist pairwise coprime linear
homogeneous polynomials
$$
l_1({\bold s}), \ldots,
l_p({\bold s}) \in {\R}\lbrack s_1, \ldots, s_I \rbrack$$
and positive integers
$k_1, \ldots, k_p$ such that  for every $j \in \{1, \ldots, p \}$
 the
linear form $l_j({\bold s})$ descends to $P_{\C}$,  $l_j({\bold s})$
does not vanish for ${\Re }({\bold s}) \in {\R}_{>0}^{I}$, and
$$
q({\bold s},{\bold t}) = f({\bold s}, {\bold t}) 
\cdot \prod_{j =1}^p l_j^{k_j}({\bold s})
$$
is analytic at ${\bold 0}$.

(iii) There exists a nonzero constant $C(f)$ and 
 a homogeneous polynomial  $q_0({\bold s})$ of degree $\mu(q)$ in 
variables $s_1, \ldots,
s_I$  such that 
$$
q({\bold s}, {\bold t}) = q_0({\bold s}) + q_1({\bold s}, {\bold t}) 
$$
and 
$$
 \frac{q_0({\bold s})}{\prod_{j =1}^p l_j^{k_j}({\bold s})} =
C(f) \cdot {\cal  X}_{\L(I)}(\psi^I({\bold s})),
$$
where 
$q_1({\bold s}, {\bold t})$ is an  analytic function at ${\bold 0}$ with
$\mu(q_1) > \mu(q_0)$,   both
functions $q_0$, $q_1$ descend  to $P_{\C}$, and
${\cal  X}_{\L(I)}$ is the ${\cal X}$-function of the cone $\L(I) =
\psi^I({\R}^I_{\geq 0}) \subset \psi^I({\R}^I)$ (see Definition 
\ref{c.func}).
}
\label{def.good}
\end{dfn}

\begin{rem}
{\rm Let $q({\bold s}, {\bold t})$ be  an arbitrary 
analytic at ${\bold 0}$ function. Collecting terms in the Taylor 
expansion of $q$, we see that   there exists a unique homogeneous
polynomial $q_0({\bold s}, {\bold t})$ and an analytic at 
${\bold 0}$ function $q_1({\bold s}, {\bold t})$ such that 
$$ q({\bold s}, {\bold t}) = q_0({\bold s}, {\bold t}) + 
q_1 ({\bold s}, {\bold t}) $$
with $\mu(q) = \mu(q_0) < \mu(q_1)$. In particular, 
the polynomial $q_0$ and the function $q_1$ in \ref{def.good} 
are uniquely defined.}
\label{unique}
\end{rem}

\begin{dfn}
{\rm If $f({\bold s},{\bold t})$ is 
good with respect to $\Gamma$ and
the set of variables $\{s_1,...,s_I\}$ as above, then the
meromorphic function
$$
\frac{q_0({\bold s})}{\prod_{j =1}^p l_j^{k_j}({\bold s})}
$$
will be called the {\em principal part of 
$f({\bold s}, {\bold t})$ at ${\bold 0}$}
and the constant  $C(f)$ the {\em principal coefficient
of $f({\bold s}, {\bold t})$ at ${\bold 0}$}. }
\end{dfn}

Suppose  that ${\rm dim }\, \psi^I ({\R}^I) \geq 2$.
Let  ${\bold \g} = ({\bold \g}^I, {\bold \g}^J) 
\in {\Z}^{I +J}$ be an element which is
not contained in $\Gamma$, $\tilde{\Gamma}: = \Gamma \oplus \Z <
{\bold \g} >$,
$\tilde{\Gamma}_{\R} := \Gamma_{\R} \oplus \R < {\bold \g} >$, 
$\tilde{P}_{\R} :=
{\R}^{I +J} /\tilde{\Gamma}_{\R}$ and $\tilde{P}_{\C} :=
{\C}^{I+J} /\tilde{\Gamma}_{\C}$. We assume that 
$\tilde{\G} \cap \R^J = {\bold 0}$ and $\tilde{\Gamma}_{\R} \cap 
{\R}_{\geq 0}^{I+J} = 0$.  We denote by 
$\tilde{\psi}$  (resp. by $\tilde{\psi}^I$)  the natural 
projection ${\C}^{I+J} \ra \tilde{P}_{\C}$ (resp. 
${\C}^{I} \ra {\C}^I/\pi^I(\tilde{\G}_{\C})$). 

The following easy statement will be helpful in the sequel:

\begin{prop}
Let $f({\bold s}, {\bold t})$ be an analytic at ${\bold 0}$ function,
$l({\bold s}) \in \R[{\bold s}]$ 
a homogeneous linear function such that $l({\bold \g}^I) \neq
0$. Assume
that $f({\bold s}, {\bold t})$ and 
$l({\bold s})$ descend to $P_{\C}$. Then
$$
\tilde{f}({\bold s}, {\bold t}) : = 
f \left({\bold s} - \frac{l({\bold s})}{l({\bold \g}^I)}
\cdot {\bold \g}^I, {\bold t}  - \frac{l({\bold s})}{l({\bold \g}^I)}
\cdot {\bold \g}^J\right)
$$
descends to $\tilde{P}_{\C}$.
\label{desc2}
\end{prop}

\begin{theo}
Let $f({\bold s},{\bold t})$ 
be a good function with respect to $\Gamma$ and
the set of variables $\{s_1,...,s_I\}$  as above,
$$
\Phi({\bold s}) = \prod_{j\;:\; l_j({\bold \g}^I)=0} l_j^{k_j}({\bold s})
$$
the product of those linear
forms $l_j({\bold s})$ $(j \in \{ 1, \ldots, p\})$ which vanish on
${\bold \g}^I$.
Assume that 
the following statements hold:

{\rm (i)} The integral
$$
\tilde{f}({\bold s}, {\bold t}) : = 
\frac{1}{2\pi i} 
\int_{{\Re }(z) = 0} f({\bold s} + z \cdot {\bold \g}^I, {\bold t} +
z \cdot {\bold \g}^J) dz
, \;\; z \in \C
$$
converges  absolutely and uniformly 
to a holomorphic function
on any compact  in the domain 
${\Re }({\bold s}) \in {\R}_{>0}^{I}$, 
${\Re }({\bold t}) \in {\R}_{> -\d_0}^{J}$;

{\rm (ii)}
There exists $\d > 0$  such that the integral
$$
\frac{1}{2\pi i} \int_{{\Re }(z) = \d}
\Phi({\bold s}) \cdot f({\bold s} + z \cdot {\bold \g}^I, 
{\bold t} + z \cdot {\bold \g}^J) dz
$$
converges absolutely and uniformly
in an  open neighborhood of ${\bold 0}$. Moreover, the
multiplicity of the meromorphic function
$$
\tilde{f}_{\d}({\bold s}, {\bold t}): = \frac{1}{2\pi i}  \int_{{\Re }(z) = \d}
 f({\bold s} + z \cdot \g) dz
$$
at  ${\bold 0}$ is at least $1 - {\dim }\, \tilde{\psi}^I ({\R}^I)$;

{\rm (iii)} For any ${\Re }({\bold s}) \in {\R}_{>0}^I$ and 
${\Re }({\bold t}) \in {\R}_{> -\d_0}^{J}$, one has 
$$
\lim_{\l \ra + \infty} 
\left( \sup_{0 \leq {\Re }(z) \leq \d, \, |{\Im }(z)| = \l  }
|f({\bold s}+ z\cdot {\bold \g}^I, {\bold t} + z \cdot {\bold \g}^J)|
\right) = 0. 
$$
Then $\tilde{f}({\bold s})$ is a good function with
respect to $\tilde{\Gamma}$ and
$\{s_1,...,s_I\}$, and $C({\tilde{f}}) = C(f)$.
\label{desc3}
\end{theo}

\noindent
{\em Proof.} By our assumption on $\tilde{\G}$, $\g^I \neq {\bold 0}$. 
We can assume that $l_j({\bold \g}^I) < 0$ for $j=1, \ldots, p_1$,
$l_j({\bold \g}^I) = 0$ for $j=p_1 +1, \ldots, p_2$, and
$l_j({\bold \g}^I) > 0$ for $j=p_2 +1, \ldots, p$. In particular, one has
\[ \Phi({\bold s}) = \prod_{j = p_1 + 1}^{p_2} l_j^{k_j}({\bold s}),
\] 
where $k_j$ $(j =p_1 +1 , \ldots, p_2)$ are some positive
integers. 
Denote by $z_j$ the solution of
the following linear equation in $z$:
\[
l_j({\bold s}) + z l_j({\bold \g}^I) = 0,\;\;j =1, \ldots, p_1.
\]

Let $U$ be the intersection of ${\R}^{I+J}_{>0}$ with an open
neighborhood of ${\bold 0}$ such that 
$ \Phi({\bold s}) \cdot
\tilde{f}_{\d}({\bold s}, {\bold t})$ 
is analytic for all $({\bold s}, {\bold t}) \in U$. 
By the property (i),  both functions
$\tilde{f}_{\d}({\bold s}, {\bold t})$ and $\tilde{f}({\bold s}, {\bold t})$
are analytic in $U$. Moreover, the
integral formulas for $\tilde{f}_{\d}({\bold s},{\bold t})$ 
and $\tilde{f}({\bold s}, {\bold t})$
show that the equalities
$\tilde{f}_{\d}({\bold u}+ iy \cdot {\bold \g}) =\tilde{f}_{\d}({\bold u})$
and $\tilde{f}({\bold u}+ iy \cdot {\bold \g}) =\tilde{f}({\bold u})$ hold
for any $y \in \R$ and ${\bold u},{\bold u}+ iy \cdot {\bold \g} \in U$. Therefore,
both functions
$\tilde{f}_{\d}({\bold s}, {\bold t})$ and $\tilde{f}({\bold s}, {\bold t})$
descend to $\tilde{P}_{\bold C}$ (see Remark \ref{desc}).
Using the properties  (i)-(iii),  we
can apply the residue theorem and obtain
\[  \tilde{f}({\bold s},{\bold t})  -
 \tilde{f}_{\d}({\bold s},{\bold t}) = 
\sum_{j=1}^{p_1} {\rm Res}_{z = z_j} f({\bold s} + z \cdot
{\g}^I,{\bold t} + z \cdot {\bold \g}^J)\]
for ${\bold s},{\bold t} \in U$.
We denote by $U({\bold \g})$ the open subset of $U$ which is 
defined by the inequalities
$$
\frac{l_j({\bold s})}{l_j({\bold \g}^I)} 
\neq \frac{l_{m}({\bold s})}{l_{m}({\bold \g}^I)}\;\;
\mbox{\rm for all $j \neq m$, $\;\;j,m \in \{ 1, \ldots, p\}$.}
$$
The open set $U({\bold \g})$ is non-empty, since we assume that
$g.c.d.(l_j, l_{m})=1$  for $j \neq m$.
Moreover, for $({\bold s},{\bold t}) \in U({\bold \g})$, we have \\
$
{\rm Res}_{z = z_j} f({\bold s} + z \cdot {\bold \g}^I,
{\bold t} + z \cdot {\bold \g}^J) = 
$
$$ =
\frac{1}{(k_j-1)!}
\left( \frac{\partial}{\partial z} \right)^{k_j-1}
\frac{l_{j}({\bold s} + z \cdot {\bold \g}^I)^{k_j}
q({\bold s} + z  \cdot {\bold \g}^I,
{\bold t} + z \cdot {\bold \g}^J)}{l_j^{k_j}
({\bold \g}^I) \cdot  \prod_{m =1}^p
l_{m}^{k_m}({\bold s} + z \cdot {\bold \g}^I) }|_{z = z_j},
$$
where
$$
 z_j = - \frac{l_j({\bold s})}{l_j({\bold \g}^I)}\, \;(j=1, \ldots, p_1) .
$$
Let
$$
f({\bold s},{\bold t}) \cdot \prod_{j =1}^p l_j^{k_j}({\bold s})
= {q}({\bold s},{\bold t}) = {q}_0({\bold s}) + {q}_1({\bold s},{\bold t})
$$
and 
$$
 \frac{q_0({\bold s})}{\prod_{j =1}^p l_j^{k_j}({\bold s})} =
C(f) \cdot  {\cal  X}_{\L(I)}({\psi}^I({\bold s})), 
$$
where ${q}_0({\bold s})$ is a uniquely determined homogeneous polynomial 
(see Remark \ref{unique}), 
${q}_0({\bold s},{\bold t})$ is an analytic at ${\bold 0}$ function
with  $\mu({q}) = \mu({q}_0) < \mu({q}_1)$ and
${\cal  X}_{\L(I)}({\psi}^I({\bold s}))$ is the ${\cal X}$-function of the cone $\L(I) =
\psi^I({\R}^{I}_{\geq 0})$.  We set
$$
R_0({\bold s}) : 
= \frac{q_0({\bold s})}{\prod_{j =1}^p {l}^{k_j}_j({\bold
s})}, \;\;
R_1({\bold s},{\bold t}) : 
= \frac{q_1({\bold s},{\bold t})}{\prod_{j =1}^p {l}_j^{k_j}({\bold
s})}.
$$
Then $\mu(f)=  \mu (R_0) < \mu (R_1)$. Moreover,
$\mu(R_0) = - {\dim }\, \psi^I(\R^I)$ (see Prop. \ref{merom}).
Define
$$
\tilde{R}_0({\bold s}):=  
\sum_{j=1}^{p_1} {\rm Res}_{z = z_j}
R_0({\bold s}+ z\cdot {\bold \g}^I)
$$
and
$$
\tilde{R}_1({\bold s},{\bold t}):=
\sum_{j=1}^{p_1} {\rm Res}_{z = z_j}
R_1({\bold s}+ z\cdot {\bold \g}^I,{\bold t} + z \cdot {\bold \g}^J).
$$
We claim that 
$$
\tilde{R}_0({\bold s}) =  C(f) \cdot  {\cal  X}_{\tilde{\L}(I)}
(\tilde{\psi}^I({\bold s})), 
$$
where 
${\cal  X}_{\tilde{\L}(I)}(\tilde{\psi}^I({\bold s}))$ 
is the ${\cal X}$-function of the cone $\tilde{\L}(I) =
\tilde{\psi}^I({\R}^{I}_{\geq 0})$.
Indeed, repeating for 
${\cal  X}_{\L(I)}(\psi({\bold s}))$
the same arguments as for $f({\bold s},{\bold t})$,  we obtain
$$
\frac{1}{2\pi i} \int_{{\Re }(z) = 0} 
{\cal  X}_{\L(I)}(\psi({\bold s} + z \cdot {\bold \g}^I)) dz
- \frac{1}{2\pi i} \int_{{\Re }(z) = \d} 
{\cal  X}_{\L(I)}(\psi({\bold s} + z \cdot {\bold \g}^I)) dz
$$
$$
= \sum_{j=1}^{k_1} {\rm Res}_{z = z_j}
{\cal  X}_{\L(I)}(\psi({\bold s} + z_j \cdot {\bold \g}^I)).
$$
Moving the contour of integration
${\Re }(z) = \d$ $(\d \ra + \infty)$, by residue theorem,
we obtain
$$
\int_{{\Re }(z) = \d} 
{\cal  X}_{\L(I)}(\psi({\bold s} + z \cdot {\bold \g}^I)) dz =0.
$$
On the other hand,
$$
{\cal  X}_{\tilde{\L}(I)}(\tilde{\psi}({\bold s})) =
\frac{1}{2\pi i}\int_{{\Re }(z) = 0}
{\cal  X}_{\L(I)}(\psi({\bold s} + z \cdot {\bold \g}^I)) dz
$$
(see Theorem \ref{char0}).

Consider the decomposition of $\tilde{f}$ into the sum:
$$
\tilde{f}({\bold s},{\bold t}) = 
\tilde{f}_{\d}({\bold s},{\bold t}) + \tilde{R}_0({\bold
s}) + \tilde{R}_1({\bold s},{\bold t}).
$$
By our assumption in (ii), $\mu(\tilde{f}_{\d}) \geq
1 -{\dim }\, \tilde{\psi}^I(\R^I)$. 
By Proposition \ref{mult4}, we
have
$\mu (\tilde{R}_1) \geq 1+ \mu(R_1)  \geq 2 + \mu(R_0)= 1- {\rm
dim}\, \tilde{\psi}^I(\R^I)$.
Using \ref{mult3}(iii),
we obtain that $\mu (\tilde{f}) = \mu (\tilde{R}_0) = - {\rm
dim}\, \tilde{\psi}^I(\R^I)$ and 
$\m(\tilde{f}_{\delta} + \tilde{R}_1) > \mu(\tilde{f})$.

By \ref{desc2}, the linear forms
\[ h_{m,j}({\bold s}):= l_{m}({\bold s} + z_j \cdot {\bold \g}^I)
= l_{m}({\bold s}) -
\frac{l_j({\bold s})}{l_{j}({\bold \g}^I)} l_{m}({\bold \g}^I), 
\; \; (j =1, \ldots, p_1,\; m \neq j) \]
and the analytic in the domain $U({\bold \g})$ functions
\[ {\rm Res}_{z = z_j} f({\bold s} + z \cdot {\bold \g}^I,{\bold t} + 
z \cdot {\bold \g}^J),\;\; j =1, \ldots, p_1 \]
and 
\[ 
{\rm Res}_{z = z_j} R_0({\bold s} + z \cdot {\bold \g}^I),
\;\; j =1, \ldots, p_1  \]
descend to $\tilde{P}_{\C}$.
For any $j \in \{ 1, \ldots, p_1\}$, let us denote
$$
Q_j ({\bold s}) = \prod_{m \neq j, m=1}^p h_{m,j}^{k_m}({\bold s}).
$$
It is clear that
$$
Q_j^{k_j}({\bold s}) \cdot {\rm Res}_{z = z_j} f({\bold s} + z \cdot
{\bold \g}^I, {\bold t} + z \cdot {\bold \g}^J)\;$$
and 
$$Q_j^{k_j}({\bold s}) \cdot {\rm Res}_{z = z_j} R_0({\bold s} + z
\cdot {\bold \g}^I)
$$
are analytic at ${\bold 0}$ and $\Phi({\bold s})$ divides
each $Q_j ({\bold s})$. Hence,  we obtain that \\
$
 \tilde{f}({\bold s},{\bold t}) \prod_{j=1}^{p_1} Q_j^{k_j}({\bold s})
  = $
$$ =  \left( \tilde{f}_{\d}({\bold s},{\bold t}) +
 \sum_{j=1}^{p_1} {\rm Res}_{z = z_j} f({\bold s} + z \cdot {\bold
\g}^I, {\bold t} + z \cdot {\bold \g}^J) \right)
\prod_{j=1}^{p_1} Q_j^{k_j}({\bold s})
$$
and
$$
 \tilde{R}_0({\bold s}) \prod_{j=1}^{p_1} Q_j^{k_j}({\bold s})
  =  \left(  
\sum_{j=1}^{p_1} {\rm Res}_{z = z_j} R_0({\bold s} + z \cdot {\bold
\g}^I) \right)
\prod_{j=1}^{p_1} Q_j^{k_j}({\bold s})
$$
are  analytic at ${\bold 0}$.

Let us define the set $\{ \tilde{l}_1({\bold s}), 
\ldots, \tilde{l}_{\tilde{p}}({\bold s}) \}$ as
a subset of pairwise coprime elements
in the set of homogeneous linear forms 
$\{ h_{m,j}({\bold s}) \}$ $(m
\in \{1, \ldots, p\}, \; j \in \{1, \ldots, p_1\})$  such that
there exist positive integers $n_1, \ldots, n_{\tilde{p}}$ and a
representation of the meromorphic functions $\tilde{f}({\bold s},{\bold t})$
and  $\tilde{R}_0({\bold s})$
as quotients
\[ \tilde{f}({\bold s},{\bold t}) =
\frac{\tilde{q}({\bold s},{\bold t})}{\prod_{j =1}^{\tilde{p}} 
\tilde{l}^{n_j}_j({\bold s})},\;\;
\tilde{R}_0({\bold s}) 
=
\frac{\tilde{q}_0({\bold s})}{\prod_{j =1}^{\tilde{p}} 
\tilde{l}^{n_j}_j({\bold s})},\]
where $\tilde{q}({\bold s},{\bold t})$ is analytic at ${\bold 0}$,
$\tilde{q}_0({\bold s})$ is a homogeneous polynomial,  and
none of the forms
$\tilde{l}_1({\bold s}), \ldots,
\tilde{l}_q({\bold s})$ vanishes for $({\bold s},{\bold t}) 
\in {\R}_{>0}^{I +J}$
(the last property can be achieved, because both functions
$\tilde{f}({\bold s},{\bold t})$ and 
$\tilde{R}_0({\bold s},{\bold t})$ are  analytic in $U$ and 
the closure of $U$ is equal to ${\R}_{\geq 0}^{I +J}$).

\noindent
Define
$$
\tilde{q}_1({\bold s},{\bold t}) = \left( \tilde{f}_{\d}({\bold s},{\bold t}) +
\tilde{R}_1({\bold s},{\bold t}) \right)
\cdot \prod_{j =1}^{\tilde{p}} \tilde{l}^{n_j}_j({\bold s}).
$$
Then
$$
\tilde{q}({\bold s},{\bold t}) = \tilde{q}_0({\bold s}) +
\tilde{q}_1({\bold s},{\bold t})
$$
where $\tilde{q}_0({\bold s},{\bold t})$ is a homogeneous polynomial
and $\tilde{q}_1({\bold s},{\bold t})$ is an analytic at ${\bold 0}$ function
such that  $\mu(\tilde{q}) = \mu(\tilde{q}_0) < \mu(\tilde{q}_1)$. 
Moreover, 
$$
 \frac{\tilde{q}_0({\bold s})}{\prod_{j =1}^{\tilde{p}}
\tilde{l}_j^{n_j}({\bold s})} =
C(f) \cdot
{\cal  X}_{\tilde{\L}(I)}(\tilde{\psi}^I({\bold s})),
$$
i.e., $\tilde{f}$ is good. 
\hfill $\Box $

\begin{dfn}
{\rm For any finite dimensional Banach space $V$ over $\R$ we 
denote by $\| \cdot\|$ a representative 
in the class of equivalent norms on $V$. For 
${\bold y}=(y_1,...,y_r)\in \R^r$ we will set  
$$
\|{\bold y}\|: =\sum_{j=1}^r |y_j|.
$$}
\end{dfn}

\noindent
The following lemma is elementary:

\begin{lem} Let $V=V_1\oplus V_2$ 
be a direct sum of finite dimensional vector spaces over $\R$,
$r_2$ is the dimension of $V_2$, and $r_2>0$.
Let $f({\bold x})$ be a complex valued
function on $V$ satisfying the inequality 
$$
|f({\bold x})|\le  \frac{c}{(1+\|{\bold x}\|)^{r_1+r_2+2\e}},
$$
for any ${\bold x}=({\bold x}_1,{\bold x}_2)\in V$ and some constants $c,\e >0$.
Let $W\subset V_2$ be a locally closed subgroup
such that $V_2/W$ is compact.
Choose any Haar measure ${\bold dw}$ on $W$. 
Then there exists a constant $c'>0$ such
that we have the estimate
$$
\int_W |f({\bold x}_1+{\bold w})|{\bold dw}
\le \frac{c'}{(1+\|{\bold x}_1\|)^{r_1+\e}}
$$
for any ${\bold x}_1 \in V_1$. 
\label{trivial}
\end{lem}

\begin{theo} 
Let $f({\bold s}, {\bold t})$ be an analytic function
for ${\Re }({\bold s})\in \R^{I}_{>0}$,
${\Re }({\bold t}) \in \R^{J}_{> -\d_0}$ $($for some $\d_0>0)$, 
$\G\subset \Z^{I+J}$ a sublattice of rank $t<I$ with 
$\G_{\R} \cap \R^J=0$ and $\G_{\R} \cap \R^{I+J}_{\geq 0} = 0$. 
Assume that there exist constants $\e,\e_0>0$ 
such that the following holds:

{\rm (i)} The function
$$
g({\bold s}, {\bold t})=s_1\cdots s_If({\bold s}, {\bold t})
$$
is holomorphic in the domain ${\Re }({\bold s}) \in \R^{I}_{> -\e}$, 
${\Re }({\bold t}) \in \R^{J}_{> -\d_0}$ 
and $C(f):=g({\bold 0})\neq 0$;

{\rm (ii)} For all $\e_1$ $($with  $0 <\e_1<\e$$)$
there exist a constant $C(\e_1)>0$ 
and an estimate 
$$
|f({\bold s}+i{\bold y}_I, {\bold t} + i{\bold y}_J) | 
\leq \frac{C(\e_1)}{(1 + \|{\bold y }\|)^{t +\e_0}}, $$
$$\; {\bold y} =( {\bold y}_I, {\bold y}_J), \;\; \|{\bold y}\| = 
\|{\bold y}_I \| + \| {\bold y}_J \|,
$$
which holds for all ${\bold s}$ such that one of the two inequalities 
$-\e < {\Re }(s_j) < \e_1$ or  ${\Re }(s_j) > \e_1$
is satisfied for every $j=1,...,I$.

Then the integral 
$$
\frac{1}{(2\pi )^t}
\int_{\G_{\R}}f({\bold s}+i{\bold y}_I, {\bold t} + i {\bold y}_J){\bold dy}
$$
is a good function with respect to $\G$ and
the set of variables $\{s_1,...,s_I\}$, and   $C(f)$ is its  
principal coefficient.   
\label{integral}
\end{theo}

\noindent
{\em Proof.} 
Without loss of generality we can assume that 
$\G$ is not contained in any of $I$ coordinate hyperplanes 
$s_j =0$ $(j =1,...,I)$, otherwise we
reduce the problem to a smaller value of $I$.  
Therefore, we can choose a basis 
${\g}^1,...,{\g}^t$ of $\G$ such that  all first $I$ coordinates of 
${\g}^u= ({\g}_I^u, {\g}_J^u) \in \Z^{I+J}$ are not
equal to $0$ for every $u=1,...,t$. 

For any non-negative integer $u\le t$ we define 
a subgroup $\G^{(u)} \subset \G$ of rank $u$ as follows:
$$
\G^{(0)}= 0;\;\;\G^{(u)}:= \bigoplus_{j=1}^u \Z<{\g}^u>,\;u=1,...,t.
$$
We introduce some auxiliary functions
$$
f^{(0)}({\bold s}, {\bold t})=f({\bold s}, {\bold t}); 
$$ $$
f^{(u)}({\bold s}, {\bold t})=
\frac{1}{(2\pi )^u}
\int_{\G^{(u)}_{\R}} 
f({\bold s}+i{\bold y}^{(u)}_I, {\bold t} +i {\bold y}^{(u)}_J)
{\bold dy}^{(u)},\;\;u=1,...,t,
$$
where ${\bold dy}^{(u)}$ is the Lebesgue measure on $\G^{(u)}_{\R}$
normalised by the lattice $\G^{(u)}$. Denote by 
$P_{\C}^{(u)}= \C^r/\G^{(u)}_{\C}$. 
By the estimate in (ii), 
$f^{(u)}({\bold s}, {\bold t})$ is a  holomorphic function in the domain
$({\rm  Re}({\bold s}),{\rm  Re}({\bold t}) ) 
\in \R^{I +J}_{>0}$ and descends to $P^{(u)}_{\C}$. 

We prove by induction 
that $f^{(u)}({\bold s}, {\bold t})$ is good with respect 
to $\G^{(u)}\subset \Z^{I+J}$ and  $\{s_1,...,s_I\}$.
By (i), $f^{(0)}({\bold s}, {\bold t})$ is good.
By induction assumption, we know that
$f^{(u-1)}({\bold s}, {\bold t})$ is good with respect to 
$\G^{(u-1)}$ and $\{s_1,...,s_I\}$. Moreover, we have
$$
f^{(u)}({\bold s}, {\bold t})=\frac{1}{(2\pi i )}\int_{{\Re }(z)=0} 
f^{(u-1)}({\bold s} + z \cdot {\g}^u_I, {\bold t} + z \cdot 
{\g}^u_J ) dz.
$$
Choose $\d_u>0$ in such a way that for every 
$j=1,...,I$ one of the following two inequalities is
satisfied:
$$
-\e < \d_u \g_j^u \ <-\e_1, \;\; \mbox{\rm or}\;\; 
\d_u\g_j^u>\e_1
$$
for some $0<\e_1<\e$.  
By (ii), the integral
$$
 f_{\d}^{(u)}({\bold s}, {\bold t})=
\frac{1}{(2\pi i )} 
\int_{{\Re }(z) = \d_u} f^{(u-1)}({\bold s} + z \cdot {\g}^u_I, 
{\bold t} + z \cdot {\g}^u_J) dz
$$
$$
=\frac{1}{(2\pi)^u} 
\int_{\G_{\R}^{(u)}} f({\bold s} + \d_u\g^{u}_I + 
i{\bold y}^{(u)}_I,{\bold t} + \d_u\g^{u}_J + 
i{\bold y}^{(u)}_J ) {\bold dy}^{(u)}
$$
converges absolutely and uniformly
in an  open neighborhood of ${\bold 0}$, i.e. the
multiplicity of $f_{\d}^{(u)}({\bold s}, {\bold t})$ is at least
$0\ge 1+ {\rm rk}\, \G^{(u)} -I$. Hence it is  holomorphic
at ${\bold 0}$ and satisfies assumption (ii) of \ref{desc3}.
By lemma \ref{trivial}, the property \ref{desc3} (iii) holds. 
Applying theorem \ref{desc3},  we conclude that
$f^{(u)}({\bold s}, {\bold t})$ is a good function with 
the principal coefficient $g({\bold 0})$. 
\hfill $\Box$

\section{Fourier analysis on algebraic tori}

Let $X_F$ be an algebraic variety over a number field $F$
and $E/F$ a finite extension of number fields. 
We shall denote by $X_E$ the $E$-variety obtained by
base change from $X_F$ and by $X(E)$ the set of
$E$-rational points of $X_F$. Sometimes we omit
the subscript in $X_E$ if the field is clear from 
the context. 

Let ${\bold G}_m= {\rm Spec}(F[x,x^{-1}])$ be the
multiplicative group scheme over  $F$.
A $d$-dimensional algebraic torus $T$ is a group scheme over $F$
such that over some finite field extension $E/F$ we have
$T_E \cong ({\bold G}_{m})^d$. We call the minimal  $E$ with this property
the splitting field of $T$. Denote by $G={\rm Gal }(E/F)$ the
Galois group of $E$  over $F$. For every 
$G$-module $A$,   $A^G$ stands for  the submodule of elements
fixed by $G$. 
For any field $E$ we denote by $\hat{T}_E$ the 
$G$-module ${\rm Hom}(T_E,{\bold G}_{m})$
of $E$-rational characters of $T$. If $E$ is the
splitting field of $T$,  we put $M: =\hat{T}_E$ and
$N: ={\rm Hom} (M,\Z)$ the dual $G$-module. We denote by $t$ the
rank of the lattice $M^G$. 

Let $T$ be an algebraic torus over a number field $F$. Denote by
${\Val }(F)$ the set of valuations of $F$ and 
by ${\Val }_{\infty}(F)$ the set of archimedian valuations.
Let $F_v$ be the
completion of $F$ with respect to $v\in {\Val }(F)$, 
${\cal V}$  an extension of $v$ to $E$,
$$
G_v\; :\;= {\rm Gal}(E_{\cal V}/F_v)\subset {\rm Gal }(E/F)
$$
the decomposition group at $v$,
$T(F_v)$ the group of $F_v$-rational points of $T$ and
$T({\cal O}_v)$ its maximal compact subgroup.
We have the canonical embeddings 
$$
\pi_v\; :\;T(F_v)/T({\cal O}_v)\hookrightarrow  N^{G_v}
$$ 
for all non-archimedian $v \in {\Val }(F)$ and 
$$
\pi_v\; :\;T(F_v)/T({\cal O}_v)\hookrightarrow  N^{G_v}_{\R}
$$ 
for all $v \in {\Val }_{\infty}(F)$. 
Denote by $\overline{x}_v$ the image of
$x_v\in T(F_v)$ in $N^{G_v}$  (resp. $N_{\R}^{G_v}$)
under $\pi_v$.

\begin{dfn}
{\rm  We call a valuation $v \in {\Val }(F)$ {\em good}, if 
the mapping $\pi_v$ is an isomorphism. We denote by $S$ a finite
subset in ${\Val }(F)$ containing ${\Val }_{\infty}(F)$ and 
all valuations $v \in {\Val }(F)$ which are not good. }
\end{dfn}

Let us recall some basic arithmetic properties of algebraic
tori over the ring of adeles ${\bold A}_F$. 
Define
$$
T^1({\bold A}_F)=\{ {\bold x}\in T({\bold A}_F) \, \mid \, 
\prod_{v\in {\Val }(F)} |m(x_v)|_v =1, \, \;
{\rm for}\,\; {\rm all}\,\; m\in M^G\}.
$$
Let ${\bold K}_T=\prod_{v\in {\Val }(F)} T({\cal O}_v) $ 
be the maximal compact subgroup
of $T({\bold A}_F)$.

\begin{prop}
The groups $T({\bold A}_F), T^1({\bold A}_F), T(F), {\bold K}_T$ have
the following properties:

{\rm (i)} $T({\bold A}_F)/T^1({\bold A}_F) \cong N_{\R}^G \cong \R^t$;

{\rm (ii)} $T^1({\bold A}_F)/T(F)$ is compact;

{\rm (iii)} $T^1({\bold A}_F)/T(F){\bold K}_T$ is isomorphic to 
the product a finite group
${\bold cl}(T)$, and a connected compact abelian group;

{\rm (iv)} $w(T)={\bold K}_T\cap T(F)$ is a finite abelian 
group of torsion elements
in $T(F)$.
\label{tori.adelic}
\end{prop}

Let $T({\cal O}) \subset T(F)$ be the subgroup 
of ${\cal O}_F$-integral points. Then $T({\cal O})$ contains 
$w(T)$,  and  ${\cal E}_T : = T({\cal O}_F)/w(T)$
has a canonical embedding, as a discrete subgroup, into 
the archimedian logarithmic space
$$
N_{\R,\infty}=
\bigoplus_{v\in {\Val }_{\infty}(F)}N_{\R}^{G_v}=
\bigoplus_{v\in {\Val }_{\infty}(F)}T(F_v)/T({\cal O}_v). 
$$
Moreover, the image of ${\cal E}_T$ 
in $N_{\R,\infty}$ is contained 
in the ${\R}$-subspace $N_{\R,\infty}^1$  
 defined as 
$$
N_{\R,\infty}^1 : = \{ \overline{x} \in N_{\R,\infty} | 
\sum_{v \in {\Val }_{\infty}(F)} m({\overline{x}}_v) = 
0\;\; \mbox{\rm for all $m \in M^G$} \}, 
$$
and the quotient $N_{\R,\infty}^1/ {\cal E}_T$ is compact. 

\begin{dfn}
{\rm Let $T$ be an algebraic torus over a number field $F$.
We define 
$$
{\cal H}_T:=(T({\bold A}_F)/T(F))^*
$$ 
as the group of topological  characters of $T({\bold A}_F)$  
which are trivial on $T(F)$. 
Define the group ${\cal D}_T$ as
$$
{\cal D}_T:=(T^1({\bold A}_F)/T(F))^*.
$$ 
Define the group ${\cal U}_T$ as:
$$
{\cal U}_T := (T^1({\bold A}_F)/T(F){\bold K}_T)^*.
$$  
We call the characters $\chi \in {\cal D}_T$ {\em discrete} and 
$\chi \in {\cal U}_T$ {\em unramified}.}
\end{dfn}

\noindent
Using \ref{tori.adelic} (i), we see that a choice of a  
splitting  of the exact sequence
$$
1 \ra T^1({\bold A}_F) \ra T({\bold A}_F) \ra 
T({\bold A}_F)/T^1({\bold A}_F) \ra 1
$$
defines  isomorphisms
$$
{\cal H}_T \cong M^G_{\R} \oplus {\cal D}_T,
$$
$$
N_{\R,\infty} = N_{\R}^G \oplus N_{\R,\infty}^1, 
$$
and 
$$
M_{\R,\infty} = M_{\R}^G \oplus M_{\R,\infty}^1, 
$$ 
where 
$$
M_{\R,\infty} = \bigoplus_{v \in {\Val }_{\infty}(F)} M_{\R}^{G_v}.
$$
and 
$M_{\R,\infty}^1$ is the minimal ${\R}$-subspace in 
$M_{\R,\infty}$ containing the image of ${\cal U}_T$ 
under the canonical mapping
$$ 
{\cal U}_T \ra  M_{\R,\infty}.
$$

>From now on we fix such a non-canonical splitting. 
This allows to consider ${\cal U}_T$ as a subgroup of ${\cal H}_T$.
By \ref{tori.adelic}, we have:

\begin{prop} 
There is an  exact sequence
$$
0 \ra {\bold cl}^*(T) \ra {\cal U}_T \ra {\cal M}_T \ra 0,
$$
where ${\cal M}_T$ is the image of the canonical 
projection of ${\cal U}_T$ to
$M_{\R,\infty}^1$  and 
${\bold cl}^*(T)$ is a finite abelian group 
dual to ${\bold cl}(T)$.
\label{ex.seq}
\end{prop}

We see from \ref{ex.seq} that a character
$\chi \in M_{\R}^G \oplus {\cal U}_T$ is determined
by its archimedian component which is an element in 
$M_{\R,\infty}$ up to a finite choice.
Denote by $y(\chi ) \in M_{\R}^G \oplus {\cal M}_T$ the image
of $\chi \in M_{\R}^G \oplus {\cal U}_T$ in $M_{\R,\infty}$. 

\noindent
For all valuations $v$ we choose Haar measures
$d\mu_v$ on $T(F_v)$ normalized by 
$$
\int_{T({\cal O}_v)} d\mu_v =1.
$$
We define the canonical measure on the group $T({\bold A}_F)$
$$
\omega = \prod_{v\in {\Val }(F)}d\mu_v.
$$
For archimedian valuations the Haar measure 
$d\mu_v$ is the pullback of the Lebesgue measure on 
$ N_{\R}^{G_v}$ under the logarithmic map 
$$
T(F_v)/T({\cal O}_v)\ra N_{\R}^{G_v}.
$$
Let ${\bold dx}$ be the Lebesgue measure on 
$T({\bold A}_F)/T^1({\bold A}_F)$.
There exists a unique Haar measure $\omega^1$ on $T^1({\bold A}_F)$
such that $\omega=\omega^1{\bold dx}$.
We define 
$$
b(T)=\int_{T^1({\bold A}_F)/T(F)}\omega^1.
$$
For any $L^1$-function $f$ on $T({\bold A}_F)$ 
and any topological character $\chi $
we denote by $\hat{f}(\chi ) $ its 
global Fourier transform with respect to $\omega$ and
by $\hat{f}_v(\chi_v ) $ the local Fourier transforms. 
We will use the following version of the Poisson formula:

\begin{theo}
Let ${\cal G}$ be a locally compact abelian group with
Haar measure $dg, {\cal G}_0\subset {\cal G} $ a closed
subgroup with Haar measure $dg_0$.
The factor group ${\cal G}/{\cal G}_0$ has a unique Haar measure $dx$
normalized by the condition $dg=dx\cdot dg_0$.
Let $f\,:\, {\cal G} \ra \C $ be an ${L}^1$-function on
${\cal G}$ and $\hat{f}$ its Fourier transform with respect
to $dg$. Suppose that $\hat{f}$ is also an ${L}^1$-function
on ${\cal G}_0^{\perp}$, where ${\cal G}_0^{\perp}$ is the group
of topological characters  $\chi $
which are trivial on ${\cal G}_0$.
Then
$$
\int_{{\cal G}_0} f(x)dg_0=\int_{{\cal G}_0^{\perp}}\hat{f}(\chi) d\chi,
$$
where $d\chi$ is the orthogonal Haar measure on ${\cal G}_0^{\perp}$
with respect to the Haar measure $dx$ on ${\cal G}/{\cal G}_0$.
\label{poi}
\end{theo}

\noindent
We will apply this formula with ${\cal G}= T({\bold A}_F)$,
${\cal G}_0= T(F)$, $dg=\omega $ and $dg_0 $ is the
discrete measure on $T(F)$. 
The Haar measure $d\chi$ induces the Lebesgue measure
on $M_{\R}^G$ normalized by the lattice $M^G\subset M_{\R}^G$ and 
the discrete measure on  ${\cal D}_T$.

\begin{dfn}{\rm
Let $T$ be an algebraic torus over
$F$ and $\overline{T(F)}$ the closure of $T(F)$
in $T({\bold A}_F)$ in the {\em direct product topology}. Define
the obstruction group to weak approximation as
$$
A(T)= T({\bold A}_F)/\overline{T(F)}.
$$
}
\label{WA}
\end{dfn}

\section{Geometry of toric varieties}

\begin{dfn} {\rm 
A complete regular $d$-dimensional fan $G$-invariant $\S$ is 
a finite set of convex rational polyhedral cones in $N_{\R}$ satisfying
the following conditions:

 (i) every cone $\s\in \S$ contains $0\in N_{\R}$;

 (ii) every face $\s'$ of a cone $ \s \in \S$ belongs to $\S$;

 (iii) the intersection of any two cones in $\S$ is a face of both 
 cones;

 (iv) $N_{\R}$ is the union of cones from $\S$;

(v) every cone $\s\in \S$ is generated by a part of a 
$\Z$-basis of $N$;

(vi) For any $g\in G$ and any $\s \in \S$, one has $g(\s)\in \S$.
 }
 \end{dfn}

A complete regular $d$-dimensional fan $\S$ defines a 
smooth toric variety
$X_{\S,E}$ as follows:
$$
X_{\S,E}=\bigcup_{\s\in \S} U_{\s}=\bigcup_{\s\in \S}
 {\rm Spec }(E[M\cap\check{ \s}]),
$$ 
where $\check{\s}\subset M_{\R}$ is the dual to $\s$ cone. We can see
that $T_E\subset U_{\s}$ for all $\s\in \S$ and that $U_0=T$.

\begin{theo}\cite{vosk1}
Let $\S$ be a complete regular $G$-invariant fan in $N_{\R}$. Assume
that the complete toric variety $X_{\S,E}$ defined over the splitting
field $E$ by $\S$ is projective. Then there exists a unique complete
algebraic variety $X_{\S,F}$ over $F$ such that its base extension
to $E$ is isomorphic to $X_{\S,E}$. 
\end{theo}

Denote by $\S(j)$ the subset of $j$-dimensional cones in $\S$
and by $N_{\s,\R}\subset N_{\R}$ the minimal linear subspace containing $\s$.
Let $\{e_1,...,e_n\}$ be the set of $1$-dimensional generators
of $\S$. Denote by $PL(\S)$ the lattice of piecewise linear integral
functions on $N$. By definition, a function $\p\in PL(\S)$ iff
$\p (N)\subset \Z$ and the restriction of $\p$ to every cone $\s\in \S$
is a linear function; equivalently, there exist 
elements $m_{\s}\in M$ such that
the restriction of $\p$ to $\s$ is given by 
$< \cdot ,m_{\s}> $
where $< \cdot ,\cdot> $ is induced from the pairing 
between $N$ and $M$. 
The $G$-action on $M$ (and $N$) induces a $G$-action on the free 
abelian group $PL(\S)$.
Let 
$$
\S(1)=\S_1(1)\cup ...\cup \S_r(1)$$
be the decomposition of $\S(1)$ into a union of 
$G$-orbits.  
A $G$-invariant piecewise linear 
function $\p\in PL(\S)^G $ is determined by the vector  
${\bold u}=(u_1,...,u_r)$,  where 
$u_i$  is the value of $\p$ on the generator of
some $1$-dimensional cone in the $G$-orbit $\S_i(1), (i=1,...,r)$. 
It will be convenient for us to
consider complex valued piecewise linear functions and to identify 
$\p=\p_{\bold u}\in PL(\S)_{\C}^G$ 
with its complex coordinates ${\bold u}=(u_1,...,u_r)\in PL(\S)_{\C}^G$.

\begin{theo} 
The toric variety $X_{\S}$ has the following properties:

(i) There is a representation of $X_{\S,E}$ as 
a disjoint union of split algebraic tori $T_{\s,E}$ of dimension
$\dim T_{\s,E}= d-\dim \s$.
For each $j$-dimensional cone $\s\in \S(j)$ we denote by 
$T_{\s,E}$ the kernel of a homomorphism $T_E\ra {\bold G}_{m,E}^j$ defined
by a $\Z$-basis of the sublattice $N\cap N_{\s,\R}$. 

(ii) The closures of $(d-1)$-dimensional tori 
corresponding to the $1$-dimensional cones 
$\R_{\ge 0}e_1,...,\R_{\ge 0}e_n\in \S(1)$ 
define divisors $\overline{T}_1,...,\overline{T}_n$. We can 
identify the lattices $PL(\S)=\oplus_{j=1}^n \Z [\overline{T}_j]$.

(iii) There is an exact sequence of $G$-modules
$$
0\ra M\ra PL(\S)\ra {\rm Pic}(X_{\S,E})\ra  0,
$$
moreover, we have ${\rm Pic}(X_{\S,F})={\rm Pic}(X_{\S,E})^G$;

(iv) The cone of effective divisors
$\L_{\rm eff}(X_{\S,F})\subset {\rm Pic}(X_{\S,F})_{\R}$
is generated by the
classes of $G$-invariant divisors 
$$
D_j=\sum_{\R_{\ge 0}e_i\in \S_j(1)} \overline{T}_i.
$$

(v) The class of the anticanonical divisor $-[K_{\S}]$ 
is given by the class of the  $G$-invariant divisor
$$
-[K_{\S}]=[D_1+...+D_r].
$$
\label{toric.geom}
\end{theo}

\begin{rem} We note that for toric varieties 
we have ${\rm Pic}(X_{\S,F})=NS(X_{\S,F})$.
\end{rem}

\begin{dfn}{\rm 
Let $\p\in PL(\S)_{\C}^G$ be a 
complex valued piecewise linear function.
Let $v\in {\Val }(F)$ be a non-archimedian valuation. 
Denote by $q_v$ the order of the residue field of $F_v$.
For $x_v\in T(F_v)$ we define the complex local height function
$$
H_{\S,v}(x_v, \p) = e^{\p (\overline{x}_v)\log q_v}. 
$$
Let $v$ be an archimedian valuation. The complex local 
height function is defined as
$$
H_{\S,v}(x_v, \p) = e^{\p (\overline{x}_v)}. 
$$
}
\end{dfn}

\begin{rem}{\rm
This provides a {\em piecewise smooth} metrization of line bundles
on the toric variety $X_{\S}$. 
One can show that this metrization is, 
in a sense, "canonical". Namely, 
an algebraic torus admits a morphism to itself 
(n-th power morphism), which extends to a compactification. 
Using the construction of Tate one can obtain a 
metrization on a line bundle by a limiting process. This metrization
coincides with ours. 
}
\end{rem}

\begin{dfn}
{\rm 
Let $x\in T(F)\subset X_{\S}(F)$
be a rational point. The global height function is defined by
$$
H_{\S}(x,\p)=\prod_v H_{\S,v}(x_v,\p).
$$
\label{height}
}
\end{dfn}

By the product formula, the function $H_{\S}(x,\p)$ 
as a function on $T(F)$ 
descends to the complexified Picard group 
${\rm Pic}(X_{\S})_{\C}$. Moreover, we have
the following 

\begin{prop} {\rm \cite{BaTschi1}} Let  
$X_{\S} $ be an smooth projective 
toric variety. For all $x\in T(F)\subset X_{\S}(F)$
the function $H_{\S}(t,\p)$
coincides with a classical height corresponding to some
metrization of the line bundle $L$ 
represented by a piecewise linear function $\p\in PL(\S)^G$. 
\end{prop}

Let $X_{\S}$ be a toric variety and ${H}_{\S}$  the height pairing. 
Clearly, it extends to a pairing
$T({\bold A}_F)\times PL(\S)_{\C}^G\ra \C$. 
Moreover, it is invariant under the maximal 
compact subgroup ${\bold K}_T=\prod_{v\in {\Val }(F)} T({\cal O}_v)$.
Therefore, its Fourier transform  
$\hat{H}_{\S}(\chi,-{\bold s})$ 
equals zero for characters 
$\chi\in {\cal H}_T$ 
which are non-trivial on ${\bold K}_T$.

By \ref{toric.geom}, we have an exact sequence of $\Z$-modules
$$
0\ra M^G\ra PL(\S)^G\ra {\rm Pic}(X_{\S,F})\ra H^1(G,M)\ra 0.
$$
It induces a surjective map of tori
$
a\,:\, \prod_{j=1}^r R_{F_j/F}({\bold G}_{m}) \ra T$
and a surjective homomorphism
$$
a\,:\,\prod_{hj=1}^r
{\bold G}_m({\bold A}_{F_j})/{\bold G}_m(F_j)\ra T({\bold A}_F)/T(F)
$$
Every character 
$\chi\in {\cal H}_T$ defines $r$ Hecke
characters $\chi_1,...,\chi_r$ of the groups 
${\bold G}_m({\bold A}_{F,j})/{\bold G}_m(F_j)$ 
by $\chi \circ a$.       
 It is known \cite{draxl}, that  ${\rm Coker}(a)$ is isomorphic
to the obstruction group to weak approximation $A(T)$ (see \ref{WA}).
Similarly, every local 
character $\chi_v$ defines local characters
$\chi_{1,v},...\chi_{r,v}$.
If $\chi$ is trivial on
${\bold K}_T$ then  all 
$\chi_j$ are trivial on the maximal compact subgroups in 
in ${\bold G}_m({\bold A}_{F_j})$, in other words, all $\chi_j$ are unramified.
Their local components for all valuations are given by
$$
\chi_{j,v}\,:\, {\bold G}_m(F_{j,v})/{\bold G}_m({\cal O}_{j,v})\ra \C^*
$$
$$
\chi_{j,v}(x_v)=|x_v|_v^{iy_{j,v}}
$$
for some real $y_{j,v}$. 
\medskip

In the remaining part of this section we recall some
estimates which will be used it the study of analytic 
properties of the height zeta function (see {\rm \cite{BaTschi1}}).   

Let us consider a Hecke character 
$\chi\in ({\bold G}_m({\bold A}_F)/{\bold G}_m(F))^*$ 
and the corresponding Hecke $L$-function $L(\chi,u)$. 
The following estimate can be proved using the Phragmen-Lindel\"of
principle \cite{rademacher}.

\begin{prop}
For all $\e>0$ there exists a constant $c_1(\e)$ such that 
for all unramified $\chi $ and all $u$ with
${\Re }(u)>1+\e$  
$$
|L(\chi,u)|<c_1(\e).
$$
For all $\e>0$ there exists a $ \delta>0$  such that for     
all unramified $\chi$
and all $u$ with ${\Re }(u)$ contained in any compact 
${\bold K}$ in the domain 
$ 0<| {\Re }\,(u) - 1| <\delta$
there exists a
constant $c({\bold K},\e)$ depending only on
${\bold K}$ and $\e$
such that
\[ | L(\chi,u) | \leq c({\bold K},\e)
(1 + |{\Im }(u)|+ \|y(\chi)\|)^{\e}.
\]
\label{m.estim}
\end{prop}

Let $T$ be an algebraic torus and $\chi\in {\cal U}_T$ an
unramified character. Denote by $\chi_v$ its local components and 
by $\chi_1,...,\chi_r$
the induced unramified Hecke characters of ${\bold G}_m({\bold A}_{F_j})$.

\begin{dfn}
{\rm Define 
$$
\zeta_{fin}(\chi,-{\bold u}):= 
\frac{\prod_{v\not\in 
{\Val }_{\infty}(F)}\hat{H}_{\S,v}(\chi_v,-{\bold u})}{
\prod_{j=1}^rL_{F_j}(\chi_j,u_j)},
$$
where for every field $F_j$ we denoted 
by $L_{F_j}(\chi_j,u) $ the standard Hecke $L$-function
of $F_j$.  
For any $\chi\in {\cal D}_T$ we define 
$$
\zeta_{\infty}(\chi,{\bold u}):=\zeta_{fin}(\chi,-{\bold u})\cdot
\prod_{v\in {\Val }_{\infty}(F)}
\hat{H}_{\S,v}(\chi,-{\bold u}).
$$
}
\label{defin}
\end{dfn}

\begin{prop} {\rm \cite{BaTschi1} }
For every $\d_0>0$ there exist constants $0<c_1<c_2$ such that 
for any ${\bold u}$ 
with  ${\Re }({\bold u})\in
\R^r_{>1/2+\d_0}$ and any $\chi\in {\cal U}_T$ 
we have 
$$
c_1<|\zeta_{fin}(\chi,-{\bold u})|<c_2.
$$
\end{prop}

\begin{prop}{\rm \cite{BaTschi1}}
Let $\chi\in {\cal U}_T$ be an unramified 
character and $y(\chi)$ its image in $M_{\R,\infty}$. 
For all $\d_0 >0$ there exists a constant $c(\d_0 )$ 
such that for any ${\bold u}$ in 
the domain ${\Re }({\bold u})\in \R^r_{>1/2+\d_0}$ 
we have the following estimate 
$$
|\prod_{v\in {\Val }_{\infty}(F)}
\hat{H}_{\S,v}(\chi,-{\bold u})| \le
\frac{c(\d_0)}{(1+\|y(\chi)\|+\|{\Im }({\bold u})\|)^{\rho+t+1}},
$$
where $\rho+t$ is the dimension of 
the real vector space $M_{\R,\infty}$.
\label{estimates-inf}
\end{prop}

\begin{coro} For any $\d_0>0$, there exists a constant $c(\d_0)$
such that for any $\chi\in {\cal U}_T$ and any ${\bold u}$ in 
the domain ${\Re }({\bold u})\in \R^r_{>1/2+\d_0}$ 
we have the following estimate:
$$
|\zeta_{\infty}(\chi,{\bold u})|\le
\frac{c(\d_0)}{(1+\|y(\chi)\|+\|{\Im }({\bold u})\|)^{\rho+t+1}}.
$$
\label{est-inf}
\end{coro}

\section{Analytic properties of height zeta functions}

\begin{dfn}{\rm Let $X_{\S}$ be a smooth projective
toric variety. Let $\p=\p_{\bold u}\in PL(\S)_{\C}^G$ be 
a complexified piecewise linear function. 
Let $Y\subset X_{\S}$ be a locally closed subset.
The height zeta function with respect to $Y$
is defined as
$$
Z_{\S}(Y, {\bold u})=\sum_{x\in Y(F)} H_{\S}(x,-{\bold u}).
$$
}
\end{dfn}

\noindent
Let us formulate the first main result. 

\begin{theo}{\rm \cite{BaTschi1}} The height zeta function 
$Z_{\S}(T,{\bold u})$
as a function on $PL(\S)_{\C}^G$ 
is holomorphic for ${\Re }({\bold u})\in \R_{>1}^r$.
Moreover, it descends to ${\rm Pic}(X_{\S})_{\C}$
and is holomorphic for ${\Re }({\bold u})$
contained in the open cone $\L_{\rm eff}^{\circ}(X_{\S})+[K_{\S}]$.
\label{convergence-cone}
\end{theo}

\begin{theo} (Poisson formula) {\rm \cite{BaTschi1,BaTschi2}}
For all ${\bold u}$ with ${\Re }({\bold s})\in \R^r_{>1}$
we have the following formula:
$$
Z_{\Sigma}(T,{\bold u})=\frac{1}{(2\pi )^t b(T)}
\int_{{\cal H}_T} \hat{H}_{\S}(\chi,-{\bold u})
d\chi,
$$
The integral converges
absolutely and uniformly  
to a holomorphic function in ${\bold u}$
in any compact in the domain ${\Re }({\bold u})\in \R^r_{>1}$.
\label{poiss}
\end{theo}

Let ${\cal L}$ be a line bundle on $X_{\S}$ metrized as above, 
such that its class $[L]$ 
is contained in the interior of the cone of effective divisors
$\L_{\rm  eff}(X_{\S})\subset {\rm Pic}(X_{\S})$. 
We have defined  $a(L)$ as 
$$
a(L):=\inf \{a\in \R \mid a[L]+ [K_{\S}] \in \L_{\rm eff}(X_{\S}) \}.
$$
By our assumptions, we have  $a(L)>0$, since
$-[K_{\S}]\in  \L^{\circ}_{\rm eff}(X_{\S})$.
Denote by 
$\L(L)$ the face of maximal codimension of the cone 
$\L_{\rm eff}(X_{\S}) $ which contains $a(L)[L]+ [K_{\S}] $. 

Let $J(L)\subset [1,...,r]$  be the set of indices such that
$[D_j]\in \L(L)$ for $j\in J(L)$ and $I(L)=[1,...,r]\backslash J(L)$.  
We set $I=|I(L)|$ and $J=|J(L)| = r - I$. Without loss of
generality, we assume that $I(L) = \{ 1, \ldots, I \}$ and 
$J(L) = \{ I+1, \ldots, r \}$.
Since $a(L)[L]+ [K_{\S}]$ is an interior point of $ \L(L)$ 
it follows that there exists a representation
$$
a(L)[L]+ [K_{\S}]=\sum_{j\in J(L)} \l_j [D_j],
$$
where $\l_j \in \Q_{>0}$. 
Therefore,
$$
[L]= \sum_{j\in J(L)} \frac{\l_j +1}{a(L)}[D_j] 
+\sum_{i\in I(L)}\frac{1}{a(L)}[D_i].
$$ 
Fix these $\l_j$ and 
choose $\e >0$ such that $2\e < \min_{j\in J(L)}\l_j$.
We denote by $\varphi_L$ the piecewise linear function from 
$PL(\S)^G_{\R}$ such that $a(L) \varphi_L(e_i) =1 $ for $i =1, \ldots, I$ 
and $a(L) \varphi_L(e_j) = 1 + \l_j$ for $j \in J(L)$.  Here $e_i$ are 
generators of one-dimensional cones $\R_{\ge 0}e_i$ in the 
$G$-orbits $\S_i(1)$.

We introduce the lattice 
$$
M_J=\{m\in M\,| \, <e,m>=0\, \;{\it for}\,\; 
\R_{\ge 0}e\in \cup_{i = 1}^I \S_i(1) \}.
$$  
Define $M_I \cong M/M_J$. The following diagram is commutative
$$
\begin{array}{ccccc}
0\ra & M_J     &\ra & \bigoplus_{j\in J(L)}\Z[G_j]   \\
     & \da     &    & \da                         \\
0\ra & M       &\ra & \bigoplus_{i=1}^r\Z[G_i]   \\
     & \da     &    & \da                         \\
0\ra & M_I     &\ra & \bigoplus_{i\in I(L)}\Z[G_i].
\end{array}
$$
The  exact sequence of $G$-modules
$$
0\ra M_J\ra M\ra M_I\ra 0
$$
induces the  exact sequence of algebraic tori 
$$
1\ra T_I\ra T\ra T_J\ra 1.
$$
It will be convenient to 
introduce new coordinates ${\bold s}=(s_i)_{i\in I(L)},
{\bold t}=(t_j)_{j\in J(L)}$ on $PL({\S})_{\C}^G$, where 
$s_i = u_i -1$ $( i = 1, \ldots, I)$, $t_j = u_{I +j} - 1 + \e$ $(j
=1, \ldots, J)$. We shall write 
$({\bold s},{\bold t})=(s_1,...,s_I, t_1,...,t_J)$.

\begin{theo}  The height zeta function 
$Z_{\S}(T,{\bold s},{\bold t})$
is good with respect to the lattice $M^G_I$ and variables 
$\{s_1,...,s_I\}$ in the domain ${\Re }({\bold s}) \in
\R_{>0}^I$, ${\Re }({\bold t}) \in
\R_{> -\d_0}^J$ for some positive $\d_0 < \e$.  
\label{Z-analytic}
\end{theo}

{\it Proof.} 
Recall that $Z_{\S}(T,{\bold u}) $ has the following integral
representation  in the domain
${\Re }({\bold u})\in \R^r_{>1}$
(\ref{poiss}):
$$
Z_{\Sigma}(T,{\bold u})=\frac{1}{(2\pi )^t b(T)}
\int_{{\cal H}_T}\hat{H}_{\S}(\chi, -{\bold u})d\chi.
$$
\noindent
Using the explicit computation of the Fourier transform of
local height functions and the absolute convergence
of the integral in the domain ${\Re }({\bold u})\in \R^r_{>1}$, 
we have 
$$
Z_{\Sigma}(T,{\bold u})=\frac{1}{(2\pi )^t b(T)}
\int_{M_{\R}^G}{\bold dy} \left(
\sum_{\chi\in{\cal U}_T}
\hat{H}_{\S}(\chi, -{\bold u}+i{\bold y}) \right),
$$
because the  local height functions are invariant under the maximal 
compact subgroups $T({\cal O}_v) \subset T(F_v)$ and   
$\hat{H}_{\S}(\chi, -{\bold u}) = 0$ for all $\chi$ which are not 
trivial on the maximal compact subgroup ${\bold K}_T$. 
By \ref{defin}, we have:
$$
\hat{H}_{\S}(\chi,-{\bold u})
=\prod_{j=1}^{r} L(\chi_j,u_j)\times \zeta_{\infty}(\chi,{\bold u}),
$$
where $\chi_1,...,\chi_r$ are unramified Hecke characters of
${\bold G}_m^1({\bold A}_{F_j})$ induced from a character 
$\chi\in {\cal U}_T$,  and
$\zeta_{\infty}(\chi,{\bold u})$ 
is a function in ${\bold u}$
which is holomorphic in the domain 
${\Re }({\bold u})\in \R^r_{>1/2+\d_0}$ (for all $\d_0>0$).

\noindent
We have
$$
{Z}_{\S}(T,{\bold s}, {\bold t})=\frac{1}{(2\pi )^t b(T)}
\int_{M_{I,\R}^G}
f_{\S}({\bold s}-i{\bold y}_I,{\bold t}-i{\bold y}_J){\bold dy}^I,
$$
where 
$$
f_{\S}({\bold s},{\bold t}):=
\sum_{\chi\in {\cal U}_T}
\prod_{i\in I(L)}L(\chi_i,s_i+1)\times b_{\S}(\chi,{\bold s},{\bold t})
$$
$$
b_{\S}(\chi,{\bold s},{\bold t})= 
$$
$$
=\int_{M_{J,\R}^G}\prod_{j\in J(L)}
L_{F_j}(\chi_j,t_j+1+ \e -iy_j) \times
\zeta_{\infty}(\chi,{\bold s}-i{\bold y}_I,{\bold t}-i{\bold y}_J)){\bold dy}^J, 
$$
${\bold dy}^I$  is the Lebesgue measure on $M_{I,\R}^G$ and  
${\bold dy}^J $ the Lebesgue measure  on $M_{J,\R}^G$.

Using the estimates \ref{trivial}, 
\ref{estimates-inf}, \ref{m.estim}, \ref{est-inf}, 
we see that the sums and integrals above
converge absolutely and uniformly to
an analytic function
in any compact in the domain ${\Re }({\bold s})\in \R^I_{>0}$
and ${\Re }({\bold t})\in \R^J_{> -\d_0}$ for some $\d_0>0$ 
$(\d_0< \e)$. 
Now the fact that 
the function $Z_{\S}(T,{\bold s},{\bold t})
$
is good with respect to 
the lattice $M_I^G\subset \Z^I$ and the variables
$(s_1,...,s_I)$ follows from \ref{integral} and the following 
statement:  

\begin{theo}
$$
\lim_{{\bold s} \ra {\bold 0}} s_1\cdots s_I f_{\S}({\bold s},{\bold 0})
$$
exists and is not equal to zero.
\label{nonzero}
\end{theo}

We divide the proof of Theorem \ref{nonzero} into a sequence 
of lemmas:

\begin{lem}
Let ${\cal U}_T(I)$ be the subgroup of ${\cal U}_T$
consisting of characters $\chi \in {\cal U}_{T}$ such that 
the corresponding Hecke characters 
$\chi_i \; \,(i =1, \ldots, I) $
are trivial. Denote 
\[ f_{\S}^I({\bold s},{\bold t})=\sum_{\chi\in{\cal U}_T(I)}
\prod_{i\in I(L)} L_{F_i}(\chi_i,s_i+1)
b_{\S}(\chi,{\bold s}, {\bold t}).\]
Then 
\[ \lim_{{\bold s} \ra {\bold 0}} s_1\cdots s_I f_{\S}({\bold s},{\bold
0}) = 
\lim_{{\bold s} \ra {\bold 0}} s_1\cdots s_I f_{\S}^I({\bold s},{\bold 0}).
\]
\end{lem}

\begin{lem} (Poisson formula)
For ${\Re }({\bold s}, {\bold t}) \in \R^I_{>0} \times \R^J_{>
-\d_0}$, 
one has:
\[ f_{\S}^I({\bold s},{\bold t}) = \int_{{\cal A}} H_{\S}(x, (-{\bold
s},- {\bold t})) d{\a},  \]
where the subgroup ${\cal A} \subset T({\bold A}_F)$ is 
defined as 
$${\cal A}:= T(F) \overline{T_I(F)}. $$
\end{lem}

{\em Proof. }
By definition of $f_{\S}^I({\bold s},{\bold t})$, we conclude  that 
this function equals to the integral of the Fourier transform 
of the adelic height function over the subgroup of 
characters $\chi$ of $T({\bold A}_F)$ which are trivial on 
$T(F)$ and such that the induced Hecke characters $\chi_i$ are trivial 
for $i \in I(L)$.  It follows from the diagram
$$
\begin{array}{cccccc}   
\prod_{i\in I(L)} {\bold G}_m({\bold A}_{F_i})/ {\bold G}_m(F_i) & \ra & 
T_I({\bold A}_F)/T_I(F)   &  \ra & A(T_I)  &\ra 0      \\
 \da   &     & \da      &      &  \da    &           \\
\prod_{i=1}^r {\bold G}_m({\bold A}_{F_i})/ {\bold G}_m(F_j)     & \ra & 
T({\bold A}_F)/T(F)       &  \ra & A(T)    &\ra 0      \\
 \da   &     & \da      &      &  \da    &           \\
\prod_{j\in J(L)} {\bold G}_m({\bold A}_{F_j})/ {\bold G}_m(F_j) & \ra & 
T_J({\bold A}_F)/T_J(F)   &  \ra & A(T_J)  &\ra 0.  
\end{array}
$$
that the common kernel of all such characters is 
$T(F) \overline{T_I(F)}$ (here we used the isomorphism 
$A(T_I) = T_I({\bold A}_F)/\overline{T_I(F)}$). 
The proof of the absolute convergence of the integral over
${\cal A}$  in the domain
${\Re }({\bold s}, {\bold t}) 
\in \R^I_{>0} \times \R^J_{> -\d_0}$
is analogous to the proof of theorem 4.2 in \cite{BaTschi2}.
 \hfill
$\Box$

\begin{lem}
The function 
\[ s_1\cdots s_If_{\S}^I({\bold s},{\bold t})  \]
extends to an analytic function in the domain 
${\Re }({\bold s}, {\bold t}) \in \R^r_{>-\d_0}$. 
\end{lem}

{\em Proof. }
The proof is similar to the proof of theorem 4.2 
in \cite{BaTschi2}.
The integral 
\[  \int_{{\cal A}} H_{\S}(x, (-{\bold s},- {\bold t})) d{\a} \]
can be estimated from above by an Euler product which is 
absolutely convergent in the domain 
${\Re }({\bold s}, {\bold t}) \in \R^r_{>-\d_0}$ times
a product of zeta functions
$\prod_{i =1}^I \zeta_{F_i}(s_i+1)$. 
\hfill $\Box$

\smallskip
\noindent
For $({\bold s},{\bold t})\in \R^r$ the function 
$H_{\S}(x, (-{\bold s},- {\bold t}))$ has values in positive
real numbers.
Therefore, to prove the non-vanishing of the constant,
it suffices to show the following:

\begin{lem}
The value of 
\[ s_1\cdots s_I \int_{\overline{T_I(F)}} H_{\S}(x, (-{\bold
s},- {\bold t})) d{\a}_I \]
at $({\bold 0}, {\bold 0})$ is positive. Here $d{\a}_I$ is the induced 
Haar measure on $\overline{T_I(F)}$
\end{lem}

{\em Proof.}
For some finite subset 
$S \subset {\Val }(F)$, we can 
split the group $\overline{T_I(F)}$ into the direct 
product
\[ \overline{T_I(F)}_S \times T_I({\bold A}_{F,S}),\]
where $\overline{T_I(F)}_S$ is the image of 
$\overline{T_I(F)}$ in the finite product 
$\prod_{v \in S} T_I(F_v)$
and 
$$T_I({\bold A}_{F,S})= T_I({\bold A}_F) \cap \prod_{v \not\in S} T_I(F_v).
$$ 
Hence,   
\[  \int_{\overline{T_I(F)}} 
H_{\S}(x, (-{\bold s},- {\bold t})) d{\a}_I = 
\]
\[
=\int_{\overline{T_I(F)}_S}  H_{\S}(x, (-{\bold
s},- {\bold t})) d{\a_S} \times \prod_{v \not\in S} 
\int_{{T_I(F_v)}}  H_{\S,v}(x, (-{\bold
s},- {\bold t})) d{\a_v}. \]
Here we denoted by $d{\a_S}$ and $d{\a_v}$ the Haar measures
induced from $d{\a}_I$.
We claim  that 
\[ \prod_{i =1}^I \zeta^{-1}_{F_i}(s_i
+1) \prod_{v \not\in S} 
\int_{{T_I(F_v)}}  H_{\S,v}(x, (-{\bold
s},- {\bold t})) d{\a_v}. \]
is an absolutely convergent  Euler product for 
${\Re }({\bold s}, {\bold t}) \in \R^r_{>-\d_0}$. 
This statement follows from the explicit calculation 
of the local integrals (see \ref{localint}).

\begin{lem}
For all good valuations $v \not\in S$, the local integral 
\[  \int_{{T_I(F_v)}}  H_{\S,v}(x, (-{\bold s},- {\bold t})) d{\a_v} = 
\prod_{i =1}^I \prod_{{\cal V}|v} \zeta_{F_i,{\cal V}}(s_i +1) 
\left( 1 + o(q_v^{-1 -\e_0}) \right) \]
for some $\e_0 > 0$ and all 
${\Re }({\bold s}, {\bold t}) \in \R^r_{>-\d_0}$
\label{localint}
\end{lem}

{\em Proof.}
Denote by $N_{\R}(I)$ the minimal ${\R}$-subspace of $N_{\R}$ spanned by
all $e$ with  $\R_{\ge 0}e$ contained in the set 
of $1$-dimensional cones in $\cup_{i\in I(L)}\S_i(1)$. 
Let $\S(L)$ be the complete $G$-invariant fan 
of cones in $N_{\R}(I)$ which consists of intersections of cones in 
$\S\subset N_{\R}$ with the subspace $N_{\R}(I)$. 
Since $\S(L)$ is not necessary 
a regular fan, we construct a new $G$-invariant fan $\tilde{\S}(L)$ 
by subdivision of cones in $\S(L)$ into regular ones using 
the method of Brylinski \cite{bryl}. 
This reduces the computation of the local intergral to 
the one made for local height functions on 
smooth toric varieties in \cite{BaTschi1}, theorem 2.2.6. 

Let $\s_1,...,\s_{\tilde{n}}$ be the set of representatives
of $G_v$-orbits in the set of
$1$-dimensional cones in ${\tilde{\S}(L)}\subset N_{\R}(I)$. 
We obtain 
$$  
\int_{{T_I(F_v)}} H_{\S,v}(x_v, (-{\bold s},- {\bold t})) d{\a_v} =
$$
$$
= Q_{\tilde{\S}(L)}(q_v^{-l_{\s_1}({\bold s}, {\bold t})},...,
q_v^{-l_{\s_{\tilde{n}}}({\bold s}, {\bold t})} )\prod_{j=1}^{\tilde{n}}
\left( 1 -  \frac{1}{q_v^{-l_{\s_j}({\bold s}, {\bold t})}} \right), 
$$
where 
$l_{\s}({\bold s}, {\bold t})$ are linear forms which are $\geq 1 -
\e_0$ in the domain   
${\Re }({\bold s}, {\bold t}) \in \R^r_{>-\d_0}$, and 
$ Q_{\tilde{\S}(L)}({\bold z}) $
is a polynomial in the variables ${\bold z}=(z_1,...,z_{\tilde{n}})$  
such that all monomials in $Q_{\tilde{\S}(L)}({\bold z}) -1$ 
have degree $\geq 2$. 
Now we notice that $l_{\s}({\bold 0}, {\bold 0}) =1$ iff $\s$ is a 
$1$-dimensional cone in $\S$ and therefore, the cone
$\R_{\ge 0}e_i$ for some $i\in I(L)$ is contained in the
$G_v$-orbit of $\s$ (see \ref{sigma}).  
\hfill 
$\Box$

\begin{lem}
The set of lattice vectors $e \in N$ such that 
$a(L)\varphi_L (e) =1$ coincides with the set 
of lattice vectors $e_i\in N_{\R}(I)$ with $\R_{\geq 0}e_i \in \S(1)$  
and $a(L)\varphi_L (e_i) =1$. 
\label{sigma}
\end{lem}

{\em Proof. } Let $e$ be a lattice point in $N$. Since $\S$ is 
complete, there exists a $d$-dimensional cone $\sigma \in \S$ such that 
$e \in \sigma$. We claim that the property $a(L)\varphi_L (e) =1$ 
implies that $e$ is a generator of a $1$-dimensional face of $\sigma$. 
Indeed, we have $a(L)\varphi_L(x) \geq \varphi_{\S}(x)$ for 
all $x \in N_{\R}$. On the other hand, 
$\sigma$ is generated by a basis of 
$N$ and $\varphi_{\S}$ has value $1$ on these generators. Hence, 
$e$ must be one of the  generators of $\sigma$.   

It remains to show that the property $a(L)\varphi_L (e_i) =1$ for 
some generator $e_i$ of a $1$-dimensional cone ${\bold R}_{\geq 0}e_i \in 
\S$ implies that $e_i \in N_{\R}(L)$. But this follows from the 
definition of $N_{\R}(L)$ as the subspace in $N_{\R}$ generated by 
all elements $e_i \in N$ such that ${\bold R}_{\geq 0}e_i \in 
\S$ and $a(L)\varphi_L (e_i) =1$. \hfill $\Box$

\begin{theo} There exists a $\d>0$ such that the zeta function
$Z_{\S}(T,{\cal L},s)$ obtained by restriction of the zeta function
$Z_{\S}(T,{\bold s})$ to the complex line 
$s[L]\in {\rm Pic}(X_{\S})_{\C}$ has a representation of the form
$$
Z_{\S}(T,{\cal L},s)= \frac{\Theta_{\cal L}(\S)}{(s-a(L))^{b(L)}}
+\frac{h(s)}{(s-a(L))^{b(L)-1}}
$$
with some function $h(s)$ which is holomorphic 
in the domain ${\Re }(s)>a(L)-\d$ and a nonzero constant
$\Theta_{\cal L}(\S)$.
\end{theo}

\section{Appendix: ${\cal X}$-functions of polyhedral cones}

Let $(A, A_{\bold R}, \L) $ be a triple consisting of
a free abelian group
$A$ of rank $k$, a $k$-dimensional real vector space
$A_{\bold R} = A \otimes {\bold R}$ containing $A$ as a sublattice of
maximal rank, and a convex $k$-dimensional 
finitely generated polyhedral cone
$\Lambda \subset A_{\R}$ such that $\Lambda \cap - \L = 0
\in A_{\R}$. Denote by  $\L^{\circ}$ the interior  of $\L$ and
by  ${\L}_{\bold C}^{\circ} = {\L}^{\circ} + iA_{\R}$
the complex  tube domain over ${\L}^{\circ}$.
Let $( A^*, A^*_{\R}, \L^*) $ be the triple
consisting of the dual abelian group
$A^* = {\rm Hom}(A, \Z)$, the dual real vector space
$A^*_{\R} = {\rm Hom}(A_{\R}, \R)$, and the  dual cone
$\L^* \subset A^*_{\R}$.
We normalize the Haar measure $ {\bold d}{\bold y}$ on $A_{\R}^*$
by the condition:
${\rm vol}(A^*_{\R}/A^*)=1$.

\begin{dfn}{\rm
 The {\em ${\cal X}$-function of} 
${\L}$ is defined as
the integral
\[  {\cal X}_{\L}({\bold s}) =
\int_{{\L}^*} e^{- \langle {\bold s}, {\bold y}
 \rangle}  {\bold d}{\bold y}, \]
where ${\bold s} \in {\L}_{\bold C}^{\circ}$.  }
\label{c.func}
\end{dfn}

\begin{prop}
One has ${\cal X}_{\L}({\bold s})$  is a rational function
$$
{\cal X}_{\L}({\bold s}) = \frac{P({\bold s})}{Q({\bold s})},
$$
where $P$ is a homogeneous polynomial,
 $Q$ is a product of all linear homogeneous forms defining
the codimension $1$ faces of
 $\L$, and ${\rm deg}\, P -
{\rm deg}\, Q = -k$.
In particular, if $(A, A_{\R},\L)=(\Z^k,\R^k,\R^k_{\ge 0})$,
then 
$$
{\cal X}_{\L}({\bold s}) = \frac{1}{s_1\cdots s_k}.
$$
\label{merom}
\end{prop}

\begin{prop} {\rm \cite{BaTschi2} }
 Let $(A, A_{\R}, \L)$ and $(\tilde{A}, \tilde{A}_{\R},
\tilde{\L})$ be two triples as above, $k = {\rm rk}\, A$ and
$\tilde{k} = {\rm rk}\, \tilde{A}$, and  $\psi\;:\; A \ra \tilde{A}$
a homomorphism of free abelian groups with a finite cokernel
${\rm Coker} (\psi )$ (i.e., the corresponding
linear mapping of real vector spaces $\psi \;:\; A_{\R} \ra
\tilde{A}_{\R}$ is surjective), and $\psi(\L) = \tilde{\L}$.
Let $\G = {\rm Ker}\, \psi \subset A$, 
${\bold d}{\bold y}$ the Haar measure
on $\G_{\R} = \G \otimes {\R}$ normalized by the condition
${\rm vol}(\G_{\R}/\G)=1$.
Then for all ${\bold s}$ with
${\Re }({\bold s}) \in \Lambda^{\circ}$
the following formula holds:
$$
{\cal  X}_{\tilde{\L}}(\psi({\bold s}))
 = \frac{1}{(2\pi)^{k-\tilde{k}}|{\rm Coker}(\psi )|}
\int_{\G_{\R}} {\cal  X}_{{\L}}
({\bold s} + i {\bold y})  {\bold dy},
$$
where $|{\rm Coker} (\psi )|$ 
is the order of the finite abelian group ${\rm Coker} (\psi )$.
\label{char0}
\end{prop}

Assume that a $\tilde{k}$-dimensional rational finite polyhedral cone
$\tilde{\Lambda} \subset \tilde{A}_{\R}$ 
contains exactly $r$ one-dimensional
faces with primitive lattice generators 
$a_1, \ldots, a_r \in \tilde{A}$.
We set $k := r$,  $A := {\Z}^r$ and 
denote by $\psi$
the natural homomorphism of lattices ${\bold Z}^r \ra \tilde{A}$
which sends the standard basis of ${\Z}^r$ into
$a_1, \ldots, a_r \in \tilde{A}$, 
so that $\tilde{\L}$ is the image
of the simplicial cone $\R^r_{\ge 0}\subset \R^r$
under the  surjective map of vector 
spaces $\psi\; : \; {\bold R}^r
\ra A_{\R}$.
Denote by $\G$ the kernel of $\psi$. 
By \ref{char0} we obtain the following:

\begin{coro}{\rm
Let ${\bold s}=(s_1,...,s_r)$ be the standard
coordinates in $\C^r$.
Then
$$
{\cal X}_{\L}(\psi({\bold s}))=
\frac{1}{(2\pi )^{r-k}|{\rm Coker}(\psi )|}
\int_{\G_{\R}}\frac{1}{\prod_{j=1,n}(s_j+iy_j)} {\bold d}{\bold y},
$$
where ${\bold dy}$ is the Haar measure on the additive
group $\G_{\R}$ normalized
by the lattice $\G$,
$y_j$ are the coordinates of ${\bold y}$ in $\R^r$, and
$|{\rm Coker} (\psi )|$  is the index of the sublattice  
in $\tilde{A}$ generated by
$a_1, \ldots, a_r$.
\label{int.formula}
}
\end{coro}

\end{document}